\documentclass[12pt, reqno]{amsart}

\usepackage{amssymb}  % AMS fonts
\usepackage{array}
\usepackage{enumitem}
\usepackage{url}
\usepackage{verbatim}
\usepackage[usenames]{color}

\theoremstyle{plain}
\newtheorem{theorem}{Theorem}[section]
\newtheorem{lemma}[theorem]{Lemma}

\theoremstyle{remark}
\newtheorem{example}[theorem]{Example}
\newtheorem{remark}[theorem]{Remark}

\numberwithin{equation}{section}

\newcommand{\seclabel}[1]{\label{sec:#1}}   % section
\newcommand{\cl}[1]{\mathrm{c\ell}(#1)}         % class

\input cyracc.def
\newfam\cyrfam
\font\tencyr=wncyr10 % cyrillic
\font\tenitcyr=wncyi10 % italic cyrillic
\def\cyr{\fam\cyrfam\twelvecyr\cyracc}
\def\itcyr{\fam\cyrfam\twelveitcyr\cyracc}

\title[Cheban loops]
{Cheban loops}

\author{J.D. Phillips}
\author{V.A. Shcherbacov}

\address[Phillips]{Department of Mathematics \& Computer Science \\
Northern Michigan University \\ Marquette, MI 49855 USA}
\email{\url{jophilli@nmu.edu}}
\urladdr{\url{http://euclid.nmu.edu/~jophilli/}}

\address[Shcherbacov]{Institute of Mathematics\\
Academy of Sciences of Moldova\\
str. Academiai 5\\
MD 2028, Chisinau, Moldova} \email{scerb@math.md}
\urladdr{http://www.math.md/imi-site/persons/Scerbacov/Scerbacov.shtml}
\urladdr{http://scerb.com}

\begin{document}

\begin{abstract}
Left Cheban loops are loops that satisfy the identity $x(xy \cdot z) = yx \cdot xz$. Right Cheban loops satisfy the mirror identity {$(z\cdot yx)x = zx\cdot xy $}. Loops that are both left and right Cheban are called Cheban loops. Cheban loops can also be characterized as those loops that satisfy the identity $x(xy \cdot z) = (y \cdot zx)x$. These loops were introduced in \cite{Ch}. Here we initiate a study of their structural properties. Left Cheban loops are left conjugacy closed. Cheban loops are weak inverse property, power associative, conjugacy closed loops; they are centrally nilpotent of class at most two.
\end{abstract}

\maketitle

\section{Introduction}
\seclabel{intro}

In \cite{Os}, Marshall Osborn showed that if $L$ is a weak inverse property loop all of whose loop isotopes also have the weak inverse property, then the nucleus of $L$ is normal and the quotient of $L$ by its nucleus is Moufang (definitions are given in the next section). He showed further that $L$ must satisfy the identity $x(yz \cdot x) = ((x \cdot yx) / x) \cdot zx$. In \cite{Ba1}, A.S. Basarab ({\cyr Basarab}) coined the term \textit{Osborn loop} to describe those loops that satisfy this identity.

The variety of Osborn loops contains as subvarieties two of the most important classes of loops---the Moufang loops and the conjugacy closed loops \cite{KKP}. One of the most important open problems in loop theory is to determine whether or not all loop isotopes of an arbitrary Osborn loop are themselves Osborn.

In \cite{Ba2}, Basarab showed that the quotient of a conjugacy closed loop by its nucleus is an abelian group. And in \cite{Ba1}, he showed that weak inverse property Osborn loops are
characterized by the identity $x(yz \cdot x) = ((1 / y) \cdot (1 / x)) \backslash 1) \cdot zx$.
He called these loops \emph{generalized Moufang loops}.

In \cite{Ch}, A.M. Cheban ({\cyr Cheban}) investigated the structure of two varieties of loops: those loops that satisfy the identity $x(xy \cdot z) = yx \cdot xz$ and those loops that satisfy the identity $x(xy \cdot z) = (y \cdot zx)x$. We call these the \emph{left Cheban identity} and the \emph{Cheban identity}, respectively. Cheban showed that  loops satisfying his second identity are generalized Moufang loops and that each of them has an abelian group as an image, which is not surprising in light of Theorem 3.9, below. He also gave an example of a loop satisfying his first identity but that was not Moufang.

As we shall see, Cheban's two varieties have many other strong structural properties; they are also intimately related to conjugacy closed loops.

Our investigations were aided by the automated reasoning tool Prover9 and the finite model builder Mace4 \cite{Mc}. We have translated most of these computer generated proofs into more ``human friendly'' form and included them in this paper. A few of the proofs, though, are quite long, and we have left them untranslated. However, each of these is posted on the first author's website \cite{Ph2}, and is clearly referenced in the proofs in this paper. We note that it is common practice to publish complicated untranslated Prover9 proofs \cite{MP, Ph3}. This is mathematically sound since the program can be made to output a simple \emph{proof object}, which can be independently verified by a short \texttt{lisp} program.

\section{Definitions}
\seclabel{definitions}

A \emph{loop} $(Q,\cdot)$ is a set $Q$ with a binary operation
$\cdot$ such that (i) for each $x\in Q$, the
\emph{left translation}
$L(x) : Q\to Q ; y\mapsto xy$
and the \emph{right translation}
$R(x) : Q\to Q ; y\mapsto yx$
are bijections, and (ii) there exists $1\in Q$ satisfying
$1\cdot x = x\cdot 1 = x$ for all $x\in Q$.
Standard references for the theory of loops are
\cite{Be, Br, Pf}.

A \emph{Moufang} loop is a loop that satisfies $xy \cdot zx = (x \cdot yz)x$. A \emph{flexible} loop satisfies $x \cdot yx = xy \cdot x$. The \emph{left alternative property}, denoted by LAP, is given by $x \cdot xy = xx \cdot y$. The RAP is the mirror identity of the LAP.

The \emph{left nucleus} of a loop $Q$ is given by $\mathrm{N}_{\lambda}(Q) = \{a : a \cdot xy = ax \cdot y, \forall x, y \in L \}$. The \emph{middle nucleus}, $\mathrm{N}_{\mu}(Q)$, and the \emph{right nucleus}, $\mathrm{N}_{\rho}(Q)$, are defined analogously. The \emph{nucleus}, then, is given by $\mathrm{N}(Q) = \mathrm{N}_{\lambda}(Q) \cap \mathrm{N}_{\mu}(Q) \cap \mathrm{N}_{\rho}(Q)$.
The \emph{commutant} of $Q$ is given by $\mathrm{C}(Q) =
\{c: \forall x \in Q,cx = xc\}.$ The \emph{center}  is the
normal subloop given by $\mathrm {Z}(Q) = \mathrm{N}(Q)\cap \mathrm{C}(Q)$.
Now, define $\mathrm {Z}_0(Q)=\{1\}$, and $\mathrm {Z}_{i+1}(Q)$, $i\ge 0$, as the preimage of $\mathrm{Z}(Q/\mathrm {Z}_i(Q))$
under the canonical projection. The loop $Q$ is \emph{(centrally) nilpotent of class} $n$,
written $\cl{Q} = n$, if $\mathrm {Z}_{n-1}(Q) < \mathrm {Z}_n(Q) = Q$.

A loop, $Q$, is \emph{left conjugacy closed}, denoted by LCC, if its left translations are closed under conjugation by left translations, i.e., if $L(x)^{-1}L(y)L(x)$ is itself a left translation for each $x,y \in Q$. This can be expressed equationally as $z\cdot yx = ((zy)/ z)\cdot zx$. \emph{Right conjugacy closed}, denoted by RCC, is the mirror identity.
A loop $Q$ is \emph{conjugacy closed}, denoted CC, if it is both LCC and RCC.
%CC-loops were introduced in 1982 by Goodaire and Robinson in \cite{GR}.
The concept of conjugacy
closedness was introduced first by Soikis \cite{SOIKIS} and later independently by Goodaire
and Robinson \cite{GR}.
In the intervening years a great deal has been discovered about their structural properties.

An especially prominent role in the analysis of CC-loops is assumed by the \emph{weak inverse property elements}, or WIP elements; i.e., those elements, $c$ such that for every $x$ in the loop we have $c(xc)^\rho = x^\rho$, where $\rho$ is the unary operation that gives the right inverse of each element $y$ in the loop, i.e., $ yy^\rho = 1$. A loop is \emph{power associative} if subloops generated by singletons are, in fact, groups. Power associative conjugacy closed loops have especially strong structural properties \cite{KK}; this variety is denoted by PACC.

A triple of bijections $(f,g,h)$ from a loop $Q_1$ to a loop $Q_2$ is called a (loop) \emph{isotopism} if $$f(x) \cdot g(y) = h(x \cdot y)$$ for every $x,y$ in $Q_1$. Note that $f$ is an \emph{isomorphism} if and only if $(f,f,f)$ is an isotopism.

Finally, an identity $\alpha=\beta$
is of \emph{Bol-Moufang type} if (i) the only operation in $\alpha$, $\beta$ is
$\cdot$, (ii) the same 3 variables appear on both sides, in the same order,
(iii) one of the variables appears twice on both sides, (iv) the remaining two
variables appear once on both sides. For instance, the Moufang law given above, $xy \cdot zx = (x \cdot yz)x$, is an identity of Bol-Moufang type. The varieties of loops classified by a single identity of Bol-Moufang type were classified in \cite{PV}. The Cheban identities are not identities of Bol-Moufang type, since the variables do not appear in the same order on both sides of the equal sign. But they do satisfy the other conditions of the definition. The varieties of loops classified by a single identity of this generalized Bol-Moufang type are classified in \cite{CHHK}.

\section{Theorems}
\seclabel{theorems}

\begin{theorem} \label{LCC_CHEB}
A loop, $Q$, is left Cheban if and only if it is LCC and $R(x)^2 = L(x)^2$ for all $x \in Q$.
\end{theorem}

\begin{proof}
In the left Cheban identity, $x(xy \cdot z) = yx \cdot xz$, let $z = 1$ to obtain $R(x)^2 = L(x)^2$. Using this, we obtain $(x \cdot xy) / x = (yx \cdot x) / x = yx$. This in turn yields $(x\backslash y)x = (x \cdot x(x\backslash y)) / x = (xy)/x$. Finally, this, together with the left Cheban law, gives $((xy)/x) \cdot xz = ((x \backslash y)x) \cdot xz = x(x(x\backslash y) \cdot z) = x \cdot yz$.

For the converse, first rearrange the LCC law to get

\

\hfill $((xy)/x) \backslash (x \cdot yz) = xz$ \hfill $(*)$

\

Next, use $R(x)^2 = L(x)^2$ to get $(x \cdot xy)/x = yx$. Now set $y = x \backslash z$ in this to get

\

\hfill $(xz) / x = (x \backslash z)x$ \hfill $(**)$

\

Combine $(*)$ and $(**)$ to get $((x \backslash y)x) \backslash (x \cdot yz) = xz$. Now use this to get $(xy) \backslash (y(yx \cdot z)) = ((y \backslash (yx))y) \backslash (y(yx \cdot z)) = yz$. Finally, multiply both sides of this by $xy$ to obtain $xy \cdot yz = (xy) \cdot ((xy) \backslash (y(yx \cdot z))) = y(yx \cdot z)$.
\end{proof}

\begin{remark}
The variety of WIP PACC-loops can be axiomatized, in the variety of loops, by the following two identities: $(xy \cdot x) \cdot xz =  x \cdot (yx \cdot x)z$ and $zx \cdot (x \cdot yx) = (z(x \cdot xy)) \cdot x$  \cite{Ph}. Left Cheban loops satisfy the first of these two identities, as the reader may check.
\end{remark}

\begin{theorem}
Let $Q$ be a left Cheban loop. Then $\mathrm{N}_{\lambda}(Q) = \mathrm{N}_{\mu}(Q)\, \unlhd \, Q$ and $\mathrm{C}(Q) \leq \mathrm{N}(Q)$. Moreover, if $a \in \mathrm{N}_{\lambda}(Q)$, then $a^2 \in \mathrm {Z}(Q)$.
\end{theorem}

\begin{proof}
The fact that $\mathrm{N}_{\lambda}(Q) = \mathrm{N}_{\mu}(Q)\unlhd Q$ follows from \cite{DRAPAL_04} Proposition 2.7, since any left Cheban loop is LCC-loop by Theorem \ref{LCC_CHEB}.

Let $a$ be a left nuclear element. Then $x \cdot ay = (x / a)a \cdot ay = a((a \cdot x / a)y) = (a(a \cdot x / a))y = ((x/a \cdot a)a)y = xa \cdot y$. That is, $a$ is a middle nuclear element. The converse is left to the reader.

Now, let $b$ be a commutant element. We have $x \cdot yb = x \cdot by = (x/b)b \cdot by = b(b(x/b) \cdot y) = b((x/b)b \cdot y) = b \cdot xy = xy \cdot b$. Thus, $b$ is a right nuclear element. For the proof that $b$ is also middle, and hence left, nuclear, see \cite{Ph2}.
Finally, let $a$ be a left nuclear element. Then, since $R(a)^2 = L(a)^2$, we have $x \cdot aa = xa \cdot a = a \cdot ax = aa \cdot x$. That is, $a^2$ is a commutant element.
\end{proof}

\begin{example}
Here is a left Cheban loop in which $1$ is nuclear but not central. This example is of minimal order.

\begin{displaymath}
\begin{array}{cccccccc}
0&1&2&3&4&5&6&7\\
1&0&3&2&5&4&7&6\\
2&4&0&6&1&7&3&5\\
3&5&1&7&0&6&2&4\\
4&2&6&0&7&1&5&3\\
5&3&7&1&6&0&4&2\\
6&7&4&5&2&3&0&1\\
7&6&5&4&3&2&1&0
\end{array}
\end{displaymath}
\end{example}

\begin{theorem}
Let $Q$ be a left Cheban loop. If $c$ is a WIP element, then $c^2$ is central. Moreover, for every $x \in Q$, $x^2$ is a WIP element, and hence, $x^4$ is central.
\end{theorem}

\begin{proof}
\cite{Ph2}.
\end{proof}

\begin{lemma}
A loop is Cheban if and only if it is both left and right Cheban.
\end{lemma}

By the left Cheban law and the right Cheban law, respectively, we have $x(xy \cdot z) = yx \cdot xz = (y \cdot zx)x$. For the converse, see \cite{Ph2}.

\begin{lemma}
Let $Q$ be a left Cheban loop. If $Q$ is a WIP loop or if $R(x^2) = L(x^2)$ for all $x \in Q$, then $Q$ is, in fact, a Cheban loop.
\end{lemma}

\begin{proof}
\cite{Ph2}.
\end{proof}

In preparation for the next theorem, recall that an \emph{extra loop} is a loop that satisfies the identity $x(y \cdot zx) = (xy \cdot z)x$. Extra loops are conjugacy closed; in extra loops, squares are nuclear \cite{KK2}.

\begin{lemma}
Let $Q$ be a left Cheban loop. If $Q$ is either flexible or satisfies the RAP, then $Q$ is, in fact, an extra loop (and, obviously, a Cheban loop).
\end{lemma}

\begin{proof}
\cite{Ph2}.
\end{proof}

\begin{lemma}
In a cancellative CC-groupoid, the following two conditions are equivalent:

(1) $R(x^2) = L(x^2)$ for all $x$

(2) $R(x)^2 = L(x)^2$ for all $x$.
\end{lemma}

The straightforward proof is left to the reader.

\begin{theorem}
A loop, $Q$, is Cheban if and only if it is conjugacy closed and satisfies $R(x^2) = L(x^2)$ for all $x \in Q$.
\end{theorem}

\begin{proof}
Combine Theorem 3.1, Lemma 3.6, and Lemma 3.8.
\end{proof}

\begin{theorem}
Cheban loops are WIP PACC-loops. Moreover, they are centrally nilpotent of class at most 2.
\end{theorem}

\begin{proof}
That Cheban loops are WIP PACC-loops is straightforward and left to the reader. By Basarab's theorem, in CC loops the commutant is contained in the center. And by Theorem 3.10, squares are are contained in the commutant, hence they are also contained in the center. So the factor of a Cheban loop by its center has exponent $2$ and is, hence, an abelian group, which finishes the proof of the theorem.
\end{proof}

\begin{remark}
A WIP PACC-loop of nilpotency class 2 need not be Cheban, as evidenced by any nonabelian group of odd order and nilpotency class 2 (since in this case, $R(x^2) = L(x^2)$ will not hold for all $x$).
\end{remark}

\begin{remark}
If $L$ is a left Cheban, right Cheban, or Cheban loop in which every element is either a square or an involution, then $L$ is an abelian group, as the reader may easily check.
\end{remark}

% Since the bibliography is set in smaller type; shrink the cyrillic
\def\cyr{\fam\cyrfam\tencyr\cyracc}
\def\itcyr{\fam\cyrfam\tenitcyr\cyracc}

\end{document}